\documentclass[12pt]{amsart}

\usepackage{amsfonts}
\usepackage{amscd}
\usepackage{amssymb}
\usepackage{enumerate}
\allowdisplaybreaks

\usepackage{amsbsy}
\usepackage{graphicx}
\usepackage{amsthm}
\usepackage{amsmath}
\usepackage{amsxtra}
\usepackage{mathrsfs}
\usepackage{bbm}

\newtheorem{theorem}{Theorem}[section]

\newtheorem{lem}[theorem]{Lemma}

\theoremstyle{definition}
\newtheorem{defn}[theorem]{Definitions}
\theoremstyle{remark}
\newtheorem{rem}[theorem]{Remark}
\numberwithin{equation}{section}

\newcommand{\R}{\mathbb R}
\newcommand{\T}{\mathbb T}
\newcommand{\Z}{\mathbb Z}
\newcommand{\Na}{\mathbb N}

\newcommand{\pa}{\partial }

\newcommand{\D}{\Delta}

\newcommand{\ap}{\alpha}

\newcommand{\vk}{\varkappa }

\title[Sharp Strichartz type estimate for harmonic oscillator]
{Sharp Strichartz type estimates for the Schr\"{o}dinger equation associated with harmonic oscillator}

\author[    P Jitendra K. Senapati and Pradeep B]{P Jitendra Kumar Senapati and Pradeep Boggarapu}
\address[P Jitendra K. Senapati.]{Department of Mathematics\\
   BITS Pilani K K Birla Goa Campus\\
    Zuarinagar, South Goa\\
403 726, Goa, India}
\email{p20180026@goa.bits-pilani.ac.in}

\address[Pradeep B.]{Department of Mathematics\\
   BITS Pilani K K Birla Goa Campus\\
    Zuarinagar, South Goa\\
403 726, Goa, India}
\email{pradeepb@goa.bits-pilani.ac.in}

\keywords{Harmonic Oscillator, Schr\"odinger equation, Strichartz inequality, Hermite expansion}
\subjclass[2020]{Primary: 35Q41, 42B37; Secondary: 42B35, 26D99.}

\begin{document}

\maketitle
\begin{abstract}
In this article we study the Schr\"odinger equation associated with Harmonic oscillator in the form of Strichartz type inequality. We give  simple proofs for  Strichartz type inequalities using purely the $L^2 \to L^p$ operator norm estimates of the spectral projections associated  harmonic oscillator proved in \cite{KT}. Our Strichartz type estimates are sharp in sense of regularity of initial data.
\end{abstract}

\section{Introduction and main results}
Consider the free Schr\"odinger equation 
\begin{equation}\label{Sch}
\begin{cases}i{u_t}(t, x)&=-\D u(t, x) \hspace{0.5cm}x\in\R^d, t\in \R\\
u(0, x)&=f(x) 
\end{cases}
\end{equation}
where $\D=\sum_{j=1}^d \frac{\pa^2}{\pa x_j^2}$,  Laplacian on $\R^d$. For $f\in L^2(\R^d)$, $e^{it\D}f$ is the unique solution to the initial value problem (\ref{Sch}). The associated Strichartz inequality reads,
\begin{equation}\label{Stel}
\|e^{it\D}f\|_{L^q((0, \infty), L^p_x(\R^d))}\leq C \|f\|_{L^2{(\R^d)}},
\end{equation}
where $d=1,~~ \mbox{and}~~ 2\leq p\leq\infty$ satisfy $\frac{1}{q}=\frac {d} {2}(\frac{1}{2}-\frac{1}{p}), d=2~~ \mbox{and}~~ 2\leq p\leq\infty,~~ \mbox{and}~~d\geq 3, 2\leq p <\frac{2d}{d-2}$ and these are the necessary conditions for \eqref{Stel}. This fundamental result proved by Keel and Tao \cite{KETA} is the source of inspiration for several extensions of \eqref{Sch} to general compact and non-compact manifolds.     
The initial value problem \eqref{Sch} has been made for the Schr\"odinger equation of the form $iu_t(t, x)+\D u(t, x) - V(x)u(t, x) = 0,$ for a suitable potential $V$ by several authors, see \cite{GV, JSS, KPV}. In particular, when $V(x)=|x|^2$,  the Strichartz inequalities have been studied in the literature see \cite{KETA, NR}. In this case the initial value problem \eqref{Sch} turns out to be an initial value problem for the Schr\"odinger equation associated with harmonic oscillator $H=-\D+|x|^2$:

\begin{equation}\label{Schh} 
\begin{cases}i{u_t}(t, x)+H u(t, x)=0 \hspace{0.5cm}x\in\R^d, t\in\R\\
u(0, x)=f(x). 
\end{cases}
\end{equation}
 The author in \cite{M3, NR} proved the Strichartz inequalities for the above problem. If $f\in L^2(\R^d),$ the solution of the initial value problem \eqref{Schh} is given by $u(t, x) = e^{itH}f(x).$ The Strichartz inequality in this case reads as,
\begin{equation}\label{Stedh}
\|e^{itH}f\|_{L^{2q}_t(\T,~  L^{2p}_{x} (\R^d))} \leq C \|f\|_{L^2(\R^d)},
\end{equation}
where $f \in L^2(\R^d),~~\T=(-\pi, \pi).$ If $p, q \geq 1$ satisfying  
\begin{equation}
\left( \frac{d-2}{d}\right)<\frac{1}{p}\leq 1 ~~\text{and}~~1\leq \frac{1}{q} \leq 2, $$ or $$ 0\leq \frac{1}{q} <1 ~~\text{and}~~ \frac{2}{q}+\frac{d}{p}\geq d.
\end{equation}
The key point for obtaining $\eqref{Stel}$ and $\eqref{Stedh}$ are the following estimates
$$\|e^{it\D}\|_{L^1(\R^d)\to L^\infty(\R^d)}\leq C |t|^{-\frac{d}{2}} ~~~~\mbox{and}~~~ \|e^{itH}\|_{L^1(\R^d)\to L^\infty(\R^d)}\leq C |\sin t|^{-\frac{d}{2}},$$
which express the dispersive property of the Schr\"odinger equation on $\R^d.$ Our main goal of this article is to study the inequality $\eqref{Stedh}$ by using $L^2\to L^p$ bounds of the spectral resolution of harmonic oscillator in the form of estimate analogous to $\eqref{Stel}.$ For all unexplained notations we direct the reader to Section \ref{Pre}.  Our main results are as follows:    
\begin{theorem}\label{mr}
Let $d\geq1,$ $2\leq p\leq\infty,~~\text{and}~~~2\leq q<\infty.$ Then we have the following Strichartz type inequality 
\begin{equation}\label{Ste}
\|e^{itH}f\|_{L_x^p(\R^d, L_t^q(\T))}\leq C_{p, q, s} \|f\|_{W^s(\R^d)},
\end{equation}
for all $s\geq \vk_{p, q},$ where $d\geq 2,$ $$\vk_{p, q}:= \begin{cases}\frac{-1}{2}(\frac{1}{2}-\frac{1}{p}) + (\frac{1}{2}-\frac{1}{q}), & 2\leq p<\frac{2(d+3)}{d+1}\\
\frac{-1}{6}+\frac{d}{6}(\frac{1}{2}-\frac{1}{p}) + (\frac{1}{2}-\frac{1}{q}), & \frac{2(d+3)}{d+1}< p \leq\frac{2d}{d-2}\\
\frac{-1}{2}+\frac{d}{2}(\frac{1}{2}-\frac{1}{p}) + (\frac{1}{2}-\frac{1}{q}), & \frac{2d}{d-2}\leq p \leq \infty\end{cases} $$
and while for $d=1,$  $$\vk_{p, q}:= \begin{cases}\frac{-1}{2}(\frac{1}{2}-\frac{1}{p}) + (\frac{1}{2}-\frac{1}{q}), & 2\leq p<4\\
\frac{-1}{6}+\frac{1}{6}(\frac{1}{2}-\frac{1}{p}) + (\frac{1}{2}-\frac{1}{q}), & 4< p \leq \infty. \end{cases} $$ For $q=2,$ and $s<\vk_{p, 2}$, the estimates in \eqref{Ste} does not hold for any non trivial initial data $f\in L^2(\R^d)$.  So, in this case, the regularity order $\vk_{p, 2}$ is sharp for $d\geq 1$ in the sense that \eqref{Ste} does not hold for all $s<\vk_{p,q}$.
\end{theorem}
For the case $p=\frac{2(d+3)}{d+1}$ and $d\geq 3$ we have the following.
\begin{theorem}\label{mr2}
Let $p=\frac{2(d+3)}{d+1}, 2\leq q< \infty$, and $d\geq 3$. Then we have the following Strichartz inequality
\begin{equation}\label{Ste2}
\|e^{itH}f\|_{L_x^p(\R^d, L_t^q(\T))}\leq C_{p, q, s} \|f\|_{W^s(\R^d)},
\end{equation}
for $\displaystyle s \geq -\frac{1}{2(d+3)}+\Big(\frac 1 2 -\frac 1 q \Big)$.
\end{theorem}

\begin{rem}
The estimates $\eqref {Ste}$ and $\eqref {Ste2}$ are different from those given in \cite{M3, NR} as the mixed norm on left hand sides of inequalities are different. 
\end{rem}

\begin{rem}~

\begin{enumerate}
\item B. Bongioanni and K. M. Rogers proved the estimate \eqref{Ste} in Theorem \ref{mr} in the case $q=2$ in \cite{BR}

\item D. Cardona in \cite{CD},  proved the estimate similar to \eqref{Ste} of Theorem \ref{mr} with initial data of general $L^p-$ Sobolev space associated with $H.$ When we compare those results for the case $L^2-$Sobolev space, our results are sharper and give better range of the exponent $p.$ 
\end{enumerate}  
\end{rem}
Organization of this article is as follows. We recall spectral decomposition of harmonic oscillator and Hermite expansion in Section \ref{Pre}. We also define Triebel-Lizorkin spaces associated to Hermite projections and their inclusion properties in Section \ref{Pre}. We prove the main results along with supported lemmas in Section \ref{Mrs}. Finally, we shall write $C$ or $C'$ to denote positive constant of independent of significant quantities, the meaning of which can change one occurrence to another. 
%

\section{Preliminaries}\label{Pre}
In this section we will recall about Hermite functions and spectral decomposition of harmonic oscillator. 
\subsection{Harmonic oscillator and Hermite Functions}
 
The harmonic oscillator (Hermite operator) is denoted by $H$ and is defined by $$H=-\D+|x|^2$$
 where $x\in \R^d$ and $$\D=\sum_{j=1}^d\frac{\pa^2}{\pa x_j^2}$$

The Hermite operator $H$ is a positive operator densely define on $L^2(\R^d)$. The Hermite functions are the eigenfunctions for the operator $H$. Hermite functions are defined in the following way. The one dimensional Hermite functions $h_k$ are defined by
$$h_k(x)={\frac{(-1)^k}{\sqrt{2^kk!\sqrt{\pi}}}}\left(\frac{d^k }{dx^k}e^{-x^2}\right)e^{\frac{x^2}{2}}$$
Then they form a complete orthonormal family in $L^2(\R)$. Now for each multi index $\mu=(\mu_1, \mu_2, \ldots, \mu_d),$ ($\mu_i$ is non-negative integer), we define the $d$-dimensional Hermite functions by tensor product: 
$$\Phi_\mu(x)=\prod_{i=1}^{d}h_{\mu_i}(x_i).$$
Then the function $\Phi_\mu$ are eigenfunctions of the Hermite operator $H$ with eigenvalue $(2|\mu|+d)$, where $|\mu|=\mu_1+\cdots+\mu_d$, i.e., 
$$H\Phi_{\mu}=(2|\mu|+d)\Phi_{\mu},$$ and they form a complete orthonormal system in $L^2(\R^d).$ Let $f\in L^2 (\R^{d})$, then the Hermite expansion is given by $$f=\sum_{\mu}\langle{f},{\Phi_\mu}\rangle\Phi_\mu(x)=\sum_{k=0}^{\infty}{P_k}f,$$
where $\langle{f},{\Phi_\mu}\rangle=\int_{\R^d}f(x)\Phi_\mu(x)dx$ and $\{P_k\}_{k\in\Na_0}$ is the orthogonal projection operator to the  eigenspace corresponding to the eigenvalue $(2k+d)$ which is given by $$P_k f(x)=\sum_{|\mu|=k}\langle f,\Phi_\mu\rangle\Phi_\mu(x).$$
Setting $\Phi_{k}(x,y)=\sum_{|\mu|=k}{\Phi_\mu(x)}{\Phi_\mu(y)},$ the Hermite projection may be written as
$$P_kf(x)=\int_{\mathbb{R}^{d}}{\Phi_{k}(x,y)f(y)dy}.$$
The functions $\Phi_{k}(x,y)$ can be obtained by the following generating function identity is called the Meheler's formula:
$$\sum_{k=0}^{\infty}{\omega^{k}\Phi_{k}(x,y)=\pi^{-\frac{d}{2}}(1-\omega^{2})^{-\frac{d}{2}}}e^{-\frac{1}{2}\frac{1+\omega^{2}}{1-\omega^{2}}({|x^2|+|y^2|})+\frac{2\omega}{1-\omega^{2}}{x}\cdot{y}},$$
for $|\omega|<1$, for more details (see \cite{ST}). Using the above formula it can be verify that $e^{itH}$ is an integral operator see \cite{JP} with the kernel $K_{it}(x, y)$ given by 
\begin{equation}\label{hsk}
K_{it}(x,y) = \frac{(2\pi)^{-d/2}}{(- i \sin 2t)^{\frac{d}{2}}} e^{\frac{-i}{2} \big(\cot 2t(|x|^2+|y|^2)-\frac{2x\cdot y}{\sin 2t}\big)}.
\end{equation}

The following theorem due to Koch and Tataru \cite{KETA} and Sanghyuk Lee et. al  \cite{JLR}, tells about the $L^2\to L^p$ operator norm of projections $P_k$ associated with harmonic oscillator.
\begin{theorem}\cite{KETA, JLR}\label{pe}
Let $k\in\Na_0, ~~d\geq1$ and $2\leq p\leq\infty.$ Then we have 
\begin{equation}\|P_kf\|_{L^p(\R^d)}\leq C k^{\vk_p}\|f\|_{L^2(\R^d)},
\end{equation}
where $d\geq 2$
$$\vk_{p}:= \begin{cases}\frac{-1}{2}(\frac{1}{2}-\frac{1}{p}), & 2\leq p<\frac{2(d+3)}{d+1}\\
\frac{-1}{6}+\frac{d}{6}(\frac{1}{2}-\frac{1}{p}), & \frac{2(d+3)}{d+1}< p \leq\frac{2d}{d-2}\\
\frac{-1}{2}+\frac{d}{2}(\frac{1}{2}-\frac{1}{p}), & \frac{2d}{d-2}\leq p \leq \infty\end{cases} $$
and while for $d=1,$  $$\vk_{p}:= \begin{cases}\frac{-1}{2}(\frac{1}{2}-\frac{1}{p}), & 2\leq p<4\\
\frac{-1}{6}+\frac{1}{6}(\frac{1}{2}-\frac{1}{p}), & 4< p \leq \infty \end{cases}. $$
and for $p=\frac{2(d+3)}{d+1}$  and $d\geq 3$ the $\vk_p=-\frac{1}{2(d+3)}.$ The exponent $\vk_p$ is sharp, in the sense that there is no $f\in L^2(\R^d),~~f\neq 0,$ such that 
\begin{equation}
\|P_kf\|_{L^p(\R^d)}\leq C' k^s\|f\|_{L^2(\R^d)},
\end{equation} 
for all $s<\vk_p.$
\end{theorem}

Given a bounded function $m$ on $\Na_0,$ we define $m(H)$ using spectral theorem by 
\begin{equation}
m(H)f=\sum^\infty_{k=0}m(2k+d)P_kf,
\end{equation}
where $f\in L^2(\R^d).$ The operator $m(H)$ is bounded on $L^2(\R^d)$ if and only if $m$ is bounded.
$$Dom (m(H)) := \big\{f\in L^2(\R^d): \sum^\infty_{k=0}|m(2k+d)|^2 \|P_kf\|^2_{L^2(\R^d)}<\infty\big\}.$$
There are some sufficient conditions on $m$ so that $m(H)$ can be extended to a bounded linear operator on $L^p(\R^d).$  For more details see \cite{ST}. In view of this, for $m(k)=e^{ikt}$ we have $Dom (m(H))=L^2(\R^d).$  In this case $$m(H)f=e^{itH}f=\sum^\infty_{k=0}e^{i(2k+d)t}P_kf,$$
for $f\in L^2(\R^d)$ and the kernel of $e^{itH}$ given by \eqref{hsk}.

\subsection{Triebel-Lizorkin spaces associated to $\{P_k\}_{k\in\Na_0}$}

\begin{defn}
Let us consider $0<p \leq \infty,~~ r\in\R~~\mbox{and}~~ 0<q\leq\infty.$ The Triebel-Lizorkin spaces associated to the family of projections  $\{P_k\}_{k\in\Na_0}$ and to the parameter $p, q~~\mbox{and}~~ r$ is defined by those complex measurable functions $f$ satisfying 
\begin{equation}
\|f\|_{F^r_{p, q}(\R^d)}=\|f\|_{F^r_{p, q}(\R^d, (P_k)_{k\in\Na_0})}:=\left\|\left(\sum^\infty_{k=0}k^{rq}|P_kf|^q\right)^\frac{1}{q}\right\|_{L^p(\R^d)}<\infty.
\end{equation} 

\end{defn}

For $s\in\R,~~ 1\leq p\leq\infty, ~~W^{s, p}(\R^d)$ denotes $L^p-$Sobolev spaces associated with $H,$ and defined by the norm 
$$\|f\|_{W^{s, p}}=\|H^sf\|_{L^p}:=\left(\int_{\R^d}|H^sf(x)|^pdx\right)^\frac{1}{p}.$$
We use the notation $W^s(\R^d)=W^{s, 2}(\R^d).$ It can be verified that 
$$\|f\|_{W^s(\R^d)}\asymp\left(\sum^\infty_{k=0}k^{2s}\|P_kf\|^2_{L^2(\R^d)}\right)^\frac{1}{2}.$$
We also have the following natural embedding properties for the Triebel-Lizorkin spaces associated to $\{P_k\}_{k\in\Na_0}$ which can be proved using similar analysis as in \cite{ERT}. 
\begin{align*}
&-~~ F^{r+\varepsilon}_{p, q_1}\hookrightarrow F^r_{p, q_1} \hookrightarrow F^r_{p, q_2} \hookrightarrow F^r_{p, \infty}, \varepsilon>0,~~0<p\leqslant\infty,~~ 0<q_1\leqslant q_2\leqslant\infty.\\
&-~~F^{r+\varepsilon}_{p, q_1} \hookrightarrow F^r_{p, q_2},~ \varepsilon>0,~~0<p\leqslant \infty,~~ 1\leqslant q_2<q_1<\infty.\\
&-~~F^0_{2, 2}=L^2~~~\mbox{and consequently, for every}~~ s\in \R,~~W^s=F^s_{2, 2}.
\end{align*}
More details of Sobolev spaces associated with Hermite expansion can be found in \cite{BT, BD, PX}.

 Throughout the article, $\T$ denotes $(-\pi, \pi)$. For all $1\leq p, q\leq\infty$, the mixed norms are defined as follows:
\begin{align*}
\|h(t, x)\|_{L^p_x(\R^d, L^q_t(\T))}&:= \|\|h(\cdot, x)\|_{L^q_t(\T)}\|_{L^p_x(\R^d)}<\infty\\
\mbox{and}~~~ \|h(t, x)\|_{L^q_t(\T, L^p_x(\R^d))}&:= \|\|h(t, \cdot)\|_{L^p_x(\R^d)}\|_{L^q_t(\T)}<\infty. 
\end{align*}

\section{Proof of main result}\label{Mrs}
In this section we will prove our main results. In the course of proof of main results, we will realize that the estimates of $L^2\rightarrow L^p$-operator norm of the projections $P_k f$ play an important role. Besides that we need the following lemmas. Let $X$ denotes subspace of $L^2(\R^d)$, spanned by Hermite functions.
\begin{lem}\label{l1}
For $1\leq p\leq \infty,$ we have 
\begin{equation}
\|u(t, x)\|_{L^p_x(\R^d, L^2_t(\T))} = \sqrt{2\pi}\|f\|_{F^0_{p, 2}(\R^d)},
\end{equation}
where $u$ is the solution of Schr\"odinger equation associated with Hermite operator \eqref{Schh} with initial data $f\in {F^0_{p, 2}(\R^d)}.$
\end{lem}
\begin{proof}
Since $X$ is dense in ${F^0_{p, 2}(\R^d)},$ it is enough to consider $f\in X,$ we can write the solution of \eqref{Schh} as
\begin{equation}\label{sd} 
u(t, x) = e^{itH}f(x)=\sum^\infty_{k=0}e^{i(2k+d)t}P_k f(x),
\end{equation}
so we have, 
\begin{align*}
\|u(t, x)\|^2_{L^2_t(\T)} &= \int^\pi_{-\pi}u(t, x)\overline{u(t, x)}dt\\
&= \int^\pi_{-\pi}\sum_{k, k'} e^{i(2k+d)t} e^{-i(2k'+d)t}P_k f(x) \overline{P_{k'} f(x)}dt\\
&=\sum_{k, k'}\left(\int^\pi_{-\pi}e^{2i(k-k')t}dt\right)P_k f(x) \overline{P_{k'} f(x)}
\end{align*}
Observe that if $k\neq k',$ then$$\int^\pi_{-\pi}e^{2i(k-k')t}dt=  \begin{cases}0, &\mbox{if}~~ k\neq k'\\
2\pi, &\mbox{if}~~ k=k' \end{cases},$$
 so we will get
\begin{align}\label{sn} 
\|u(t, x)\|^2_{L^2_t(\T)}&=\sum^\infty_{k=0}2\pi|P_kf|^2=2\pi \sum^\infty_{k=0}|P_kf|^2.
\end{align}
Consider 
\begin{align}\label{mdn}
\|u(t, x)\|_{L_x^p(\R^d, L^2_t(\T))}&= \|\|u(t, x)\|^2_{L^2_t(\T)}\|_{L_x^p(\R^d)}\\
&=\sqrt{2\pi}\left\|\left( \sum^\infty_{k=0}|P_kf|^2\right)^{\frac{1}{2}}\right\|_{L_x^p(\R^d)} \nonumber \\
&=\sqrt{2\pi}\|f\|_{F^0_{p, 2}(\R^d)},\nonumber
\end{align}
which completes the proof of the lemma.
\end{proof}

\begin{lem}\label{l2}
Let $0 <p\leq \infty,$ $2\leq q <\infty$ and $s_q :=\frac{1}{2}-\frac{1}{q}.$ Then 
\begin{equation}
C_q' \|f\|_{F^0_{p, 2}(\R^d)}\leq \|u(t, x)\|_{L_x^p(\R^d, L^q_t(\T))} \leq C_{q, s} \|f\|_{F^s_{p, 2}(\R^d)},
\end{equation}
holds for every $s\geq s_q.$
\end{lem}
\begin{proof}
Let us consider $f\in X.$ In order to estimate the norm $\|u(t, x)\|_{L_x^p(\R^d, L^q_t(\T))}$ we can use the Wainger Sobolev embedding Theorem:
\begin{equation}
\left\|\sum_{n\in \Z, n\neq 0}|n|^{-\ap}\widehat{F}(n)e^{int}\right\|_{L^q(\T)}\leq C \|F\|_{L^r(\T)}, ~~\ap :=\frac{1}{r}-\frac{1}{q}
\end{equation}
For $s_q :=\frac{1}{2}-\frac{1}{q}~~ \mbox{and}~~ f\in X$ we have
\begin{align*}
\|u(t, x)\|_{L^q(\T)} &= \left\|\sum^\infty_{k=0}e^{i(2k+d)t}P_k f(x)\right\|_{L^q(\T)}\\
&= \left\|\sum^\infty_{k=0}(2k+d)^{-s_q}e^{i(2k+d)t}(2k+d)^{s_q}P_k f(x)\right\|_{L^q(\T)}\\
&\leq C_q \left\|\sum^\infty_{k=0}(2k+d)^{s_q}e^{i(2k+d)t}P_k f(x)\right\|_{L^2(\T)}\\
&=C_q \left(\sum^\infty_{k=0}(2k+d)^{2s_q}|P_k f(x)|^2\right)^\frac{1}{2}\\
&\leq C_q  \left(\sum^\infty_{k=0} k^{2s_q}|P_k f(x)|^2\right)^\frac{1}{2}.
\end{align*}
Thus, we get
\begin{align*}
\|u(t, x)\|_{L^p_x(\R^d, L^q(\T)} &\leq C_q  \left\|\left(\sum^\infty_{k=0} k^{2s_q}|P_k f(x)|^2\right)^\frac{1}{2}\right\|_{L^p_x(\R^d)}\\
&=C_q \|F\|_{F^{s^q}_{p, 2}(\R^d)}\leq C_{q, s} \|F\|_{F^{s}_{p, 2}(\R^d)}
\end{align*}
Last inequality in the above is due to the embedding $F^s_{p, 2}\hookrightarrow F^{s_q}_{p, 2}$ for every $s\geqslant s_q.$  From Lemma \ref{l1} we know 
\begin{equation}\label{a1}
\|f\|_{F^0_{p, 2}(\R^d)} = C \|u(t, x)\|_{L^p_x(\R^d, L^2_t(\T))}.
\end{equation}
Since  $\|u(t, x)\|_{L^2_t(\T)}\leq C_q' \|u(t, x)\|_{L^q_t(\T)}$ for $2\leq q< \infty$, in view of \eqref{a1}, we will get
$$\|f\|_{F^0_{p, 2}(\R^d)}\leq C_q' \|u(t, x)\|_{L^p_x(\R^d, L^2_t(\T))}.$$ 
This proves our Lemma.       
\end{proof}

%

Now, with the analysis developed above and by using Theorem \ref{pe} we will provide a short proof for the main results. 
\begin{proof}[Proof of Theorem \ref{mr} and \ref{mr2}] We have that from Lemma \ref{l2},
$$C_{q}' \|f\|_{F^0_{p, 2}(\R^d)}\leq \|u(t, x)\|_{L_x^p(\R^d, L^q_t(\T))} \leq C_{q, s} \|f\|_{F^s_{p, 2}(\R^d)},$$
so, in view of the above it is enough to estimate the ${F^s_{p, 2}(\R^d)}$-norm of the initial data $f,$ in terms of $L^2$-Sobolev norm $ \|f\|_{W^s}$ of $f$,  for every $s\geq\vk_{p, q}.$ Moreover, by the embedding $W^s \hookrightarrow W^{\vk_{p, q}}$ for every $s\geq\vk_{p, q}.$ It is enough to prove that $$\|f\|_{F^s_{p, 2}(\R^d)}\leq C \|f\|_{W^{\vk_{p, q}}}.$$
Now Consider for $2\leq p<\infty$ together with the Minkowski integral inequality and Theorem \ref{pe}, we get
\begin{align*}
\|f\|_{F^s_{p, 2}(\R^d)} &= \left\|\left(\sum^\infty_{k=0}k^{2{s_q}}|P_kf|^2\right)^\frac{1}{2}\right\|_{L^p_x(\R^d)}\\
&\leq \left(\sum^\infty_{k=0}k^{2{s_q}}\|P_kf\|^2_{L^p(\R^d)}\right)^\frac{1}{2}\\
&\leq C \left(\sum^\infty_{k=0}k^{2({s_q+\vk_p})}\|P_kf\|^2_{L^2(\R^d)}\right)^\frac{1}{2}\\
&\leq C \left(\sum^\infty_{k=0}(2k+d)^{2({s_q+\vk_p})}\|P_kf\|^2_{L^2(\R^d)}\right)^\frac{1}{2}\\
&= C \|f\|_{W^{\vk_{p, q}}}, 
\end{align*}   
  where we have used that $\vk_{p, q}=s_q+\vk_p.$ Let us note that the previous estimates are valid for $p=\infty.$  We prove the sharpness results for the case $q=2$ by choosing $f=P_kg$, where $g$ is arbitrary $g\ne 0\in L^2(\R^d)$. Suppose the inequality \eqref{Ste} or \eqref{Ste2} holds for some $s<\vk_p=\vk_{p, 2}$, then for the choice of $f=P_kg$, we get that
$$\|e^{itH}f\|_{L^p_x(\R^d, L^2(\T))}=\|P_kg\|_{L^p(\R^d)}\leq \|f\|_{W^s}=k^s \|P_kg\|_{L^2(\R^d)}\leq k^s \|g\|_{L^2(|R^d)}$$
\end{proof}
which leads to the improvement of the estimates for the spectral projection operators $P_k$'s which is not possible.

\section*{Acknowledgments}
Both the authors are grateful to the Department of Mathematics, BITS Pilani K K Birla Goa Campus for the facilities utilized. The authors thank Prof. S. Thangavelu, IISc Bangalore for the fruitful discussions.


\begin{thebibliography}{10}



\bibitem{BR} B. Bongioanni and K. M. Rogers, \textit{Regularity of the Schr\"odinger equation for the harmonic oscillator}, Ark. Mat., 49 (2011), 217-238.

\bibitem{BT} B. Bongioanni and Torrea, J. L, \textit{Sobolev spaces associated with the harmonic oscillator},  Proc. Indian Acad. Sci. Math. Sci. 116 (2006), no. 3, 337–360.

 \bibitem{BD} Bui The Anh and Duong, X. T,  \textit{Besov and Triebel-Lizorkin spaces associated with Hermite operators}, J. Fourier Anal. Appl. 21 (2015), no. 2, 405-448. 

\bibitem{CD} Duv\'an Cardona, \textit{$L^p$-estimates for a Schr\"{o}dinger equation associated with the harmonic oscillator}, Electronic Journal of  Differential Equations. Vol. 2019 (2019), No. 20, pp. 1-10.

\bibitem{ERT} N. Erlan, M. Ruzhansky and S. Tikhonov, \textit{Nikolskii inequality and Besov, triebel-Lizorkin, Wiener and Beuriling spaces on compact homogeneous manifolds}, Ann. Sc. Norm. Super. Pisa Cl. Sci. (5) 16 (2016), no. 3, 981-1017. 

\bibitem{GV} J. Ginber and G. Velo, \textit{The glogal Cauchy problem for the nonlinear Sch\"odinger equation revisited,} Ann. Inst. H. Poincar\'e Anal. Non Lin\'ear 2(4), 309-327.

\bibitem{JLR} E. Jeong, S. Lee and J. Rau, \textit{End Point eigenfunction bounds for the Hermite operator}, arxiv:2205.03036v1.

\bibitem{JSS} J. L. Journe, A. Soffer and C. D. Sogge, \textit{Decay estimates for Schr\"odinger operators,} Comm. Pure Appl. Math. 44(5), 573-604 (1991).

\bibitem{KT} M. Keel and T. Tao, \textit{Endpoint Strichartz estimates}, Amer. J. Math., 120 (1998), pp. 955-980.

\bibitem{KPV} C. E. Kenig, G. Ponce and L. Vega, \textit{Smoothing effects and local exsitance theory for the generalized nonlinear Schr\"odinger equations,} Invent. Math. 134(3), 489-545 (1998).

\bibitem{KETA} H. Koch and D. Tataru, \textit{Lp eigenfunction bounds for the Hermite operator}, Duke Math.J. 128 (2005), 369–392.


\bibitem{M3} Mejjaoli, H, \textit{Dunkl-Schr\"{o}dinger Equations with and without Quadratic Potentials}. Serdica Mathematical Journal 37.2 (2011): 113-142.

\bibitem{NR}  A. K. Nandakumaran and P. K. Ratnakumar, \textit{Schr\'odinger equation and the oscillatory semigroup for the Hermite operator}, J. Funct. Anal. 224(2), 371-385(2005).

\bibitem{PX} Petrushev, P and Xu, Y, \textit{Decomposition of spaces of distributions induced by Hermite expansions}. J. Fourier Anal. Appl. 14 (2008), no. 3, 372  - 414.

\bibitem{JP} P J. K. Senapati and Pradeep B., \textit{Strichartz inequality for orthonormal functions associated with Dunkl Laplacian and Hermite-Schr\"{o}dinger operators}, arXiv:2208.13024v2. 

\bibitem{ST} S. Thangavelu, \textit{Lectures on Hermite and Laguerre Expansions}, Mathematical Notes, vol. 42. Princeton University Press, Princeton (1993).


\end{thebibliography}
\end{document}